\documentclass[journal]{IEEEtran}

\usepackage{myStyleIEEE}
\usepackage[margin=0.55in]{geometry}
\usepackage{enumitem}
\usepackage{xfrac}
\usepackage{lipsum}

\setlist[itemize]{leftmargin=*}

\makenomenclature

\IEEEoverridecommandlockouts
\usepackage{mathtools,nccmath}
\newcolumntype{P}[1]{>{\centering\arraybackslash}p{#1}}
\begin{document}

\title{A Stochastic-Robust Approach for Resilient Microgrid Investment Planning Under Static and {\color{black}Transient} Islanding Security Constraints}

\renewcommand{\theenumi}{\alph{enumi}}

\newcommand{\uros}[1]{\textcolor{magenta}{$\xrightarrow[]{\text{UM}}$ ``#1''}}
\newcommand{\agnes}[1]{\textcolor{blue}{$\xrightarrow[]{\text{AN}}$ ``#1''}}
\newcommand{\shahab}[1]{\textcolor{red}{$\xrightarrow[]{\text{SD}}$ ``#1''}}
\newcommand{\petros}[1]{\textcolor{purple}{$\xrightarrow[]{\text{PA}}$ ``#1''}}

\author{Agnes~Marjorie~Nakiganda,~\IEEEmembership{Student~Member,~IEEE,
}
Shahab~Dehghan,~\IEEEmembership{Senior~Member,~IEEE,}
Uros~Markovic,~\IEEEmembership{Member,~IEEE,}
Gabriela~Hug,~\IEEEmembership{Senior~Member,~IEEE,}
and~Petros~Aristidou,~\IEEEmembership{Senior~Member,~IEEE}\vspace{-0.5cm}
        
}

\maketitle
\IEEEpeerreviewmaketitle

\begin{abstract}
    When planning the investment in Microgrids (MGs), {\color{black} usually} static security constraints {\color{black} are included} to ensure their resilience and ability to operate in islanded mode. However, unscheduled islanding events may trigger cascading disconnections of Distributed Energy Resources (DERs) inside the MG due to the transient response, leading to a partial or full loss of load. In this paper, a min-max-min, hybrid, stochastic-robust investment planning model is proposed to obtain a resilient MG considering both High-Impact-Low-Frequency (HILF) and Low-Impact-High-Frequency (LIHF) uncertainties. The HILF uncertainty pertains to the unscheduled islanding of the MG after a disastrous event, and the LIHF uncertainty relates to correlated loads and DER generation, characterized by a set of scenarios. The MG resilience under both types of uncertainty is ensured by incorporating static and transient islanding constraints into the proposed investment model. The inclusion of transient response constraints leads to a min-max-min problem with a non-linear dynamic frequency response model that cannot be solved directly by available optimization tools. Thus, in this paper, a three-stage solution approach is proposed to find the optimal investment plan. The performance of the proposed algorithm is tested on the CIGRE 18-node distribution network. 
\end{abstract}
\vspace{-0.4cm}

\vspace{-0.3cm}

 \makenomenclature
 \renewcommand\nomgroup[1]{%
  \item[\bfseries
   \ifstrequal{#1}{A}{Functions}{%
   \ifstrequal{#1}{B}{Indices}{%
   \ifstrequal{#1}{C}{Parameters}{%
   \ifstrequal{#1}{D}{Sets}{%
   \ifstrequal{#1}{E}{Symbols}{%
   \ifstrequal{#1}{F}{Variables}{%
   \ifstrequal{#1}{G}{Vectors}{}}}}}}}%
 ]}
 \nomenclature[A, 01]{\scalebox{0.8}{$\Theta^{\mathbf{gm,opr}}$}}{Total operational costs in grid-connected mode $[\$]$.}
 \nomenclature[A, 02]{\scalebox{0.8}{$\breve\Theta^{\mathbf{im,opr}}$}}{Vector-valued function of total penalty costs of disconnecting loads from MG in islanded mode $[\$]$.}
 \nomenclature[A, 01]{\scalebox{0.8}{$\Theta^{\mathbf{im,opr}}_{to}$}}{Total penalty costs of disconnecting loads from MG at hour $t$ in representative day $o$ in islanded mode $[\$]$.}
 \nomenclature[A, 03]{\scalebox{0.8}{$\Theta^{\mathbf{inv}}$}}{Total investment costs $[\$]$.}

 \nomenclature[B, 01]{$g$}{Index of generators, $g\in\{c,d,i,v\}$.}
 \nomenclature[B, 02]{$n$}{Index of nodes, $n''$/$n'$ {\color{black}being} a node before/after node $n$.}
 \nomenclature[B, 03]{$o$}{Index of representative days.}
 \nomenclature[B, 04]{$t$}{Index of hours.}
  \nomenclature[B, 05]{$\psi$}{Index of iterations.}

 \nomenclature[C, 00]{$\color{black}c_{go}$}{Daily capacity factor of generator $g$ in representative day $o$.}
 \nomenclature[C, 01]{$D$}{Normalized damping constant of all generators $[\mathrm{p.u.}]$.} 
 \nomenclature[C, 02]{$D_g$}{Normalized damping constant of SGs $[\mathrm{p.u.}]$.} 
  \nomenclature[C, 03]{$D_{v}$}{Virtual damping constant of VSM-based CIG $v$ $[\mathrm{p.u.}]$.}
 \nomenclature[C, 04]{$d_{nto}^\mathrm{p_c/q_c}$}{Constant part of active/reactive load $d_{nto}^\mathrm{p/q}$ $[\mathrm{kWh}]$.}
 \nomenclature[C, 06]{$e_{no}$}{Flexible energy demand of node $n$ in representative day $o$ $[\mathrm{kWh}]$.}
 \nomenclature[C, 07]{$e_{to}^\mathrm{b/s}$}{Buying/selling price of electricity from/to the main grid at hour $t$ in representative day $o$ $[\$/\mathrm{kWh}]$.}
 \nomenclature[C, 08]{$F_{i}$}{Fraction of the total power generated by the turbine of SG $i$, $F_{g}$ {\color{black}being} the weighted average of all SGs $[\mathrm{p.u.}]$.}
 \nomenclature[C, 09]{$fc_{n}$}{Penalty cost of shifting demand at node $n$ $[\$/\mathrm{kWh}]$.}
 \nomenclature[C, 10]{$ic_g$}{Annualized investment cost of generator $g$ $[\$]$.}
 \nomenclature[C, 11]{$ic_{nn'}$}{Annualized investment/reinforcement cost of a line connecting nodes $(n,n')$ $[\$]$.}
  \nomenclature[C, 12]{$K_d$}{Power gain factor of droop-based CIG $d$ $[\mathrm{p.u.}]$.}
  \nomenclature[C, 13]{$K_i$}{Mechanical power gain factor of SG $i$ $[\mathrm{p.u.}]$.}
  \nomenclature[C, 14]{$M$}{Normalized inertia constant of all SGs and CIGs $[\mathrm{s}]$.}
  \nomenclature[C, 15]{$M_{g}$}{Normalized inertia constant for the CoI of SGs $[\mathrm{s}]$.}
  \nomenclature[C, 16]{$M_{v}$}{Virtual inertia constant of CIG $v$ with VSM control $[\mathrm{s}]$.}
 \nomenclature[C, 17]{$mc_{g}$}{Marginal cost of generator $g$ $[\$/\mathrm{kWh}]$.}
 \nomenclature[C, 18]{$pc_{n}$}{Penalty cost of disconnecting demand at node $n$ $[\$/\mathrm{kWh}]$.}

 \nomenclature[C, 19]{$R_{d}$}{Droop of CIG $d$ with droop control $[\mathrm{\%}]$.}
 \nomenclature[C, 20]{$R_i$}{Droop of SG $i$, $R_{g}$ {\color{black}being} the weighted average of all SGs $[\mathrm{\%}]$.}
 \nomenclature[C, 21]{$r_g^\mathrm{d/u}$}{Ramp-down/ramp-up limit of generator $g$ $[\mathrm{kW}/\mathrm{h}]$.}
 \nomenclature[C, 22]{$r_{n^{\prime}n}$}{Resistance of a line connecting nodes $(n,n')$ $[\mathrm{\Omega}]$.}
 \nomenclature[C, 23]{${\overline s}_{n^{\prime}n}$}{Capacity of a line connecting nodes $(n,n')$ $[\mathrm{kVA}]$.}
 \nomenclature[C, 24]{$T_{{d}/{v}}$}{Time constant of CIG $d/v$ with droop/VSM control $[\mathrm{ms}]$.}
 \nomenclature[C, 25]{$T_i$}{Turbine time constant of SG $i$ $[\mathrm{s}]$.}
 \nomenclature[C, 26]{$x_{n^{\prime}n}$}{Reactance of a line connecting nodes $(n,n')$ $[\mathrm{\Omega}]$.}
 \nomenclature[C, 27]{$z_{nn'}^0$}{Initial status of a line connecting nodes $(n,n')$ (i.e., $1/0$: built/not-built).}
 \nomenclature[C, 28]{$\alpha$}{Scaling factor.}
 \nomenclature[C, 29]{$\zeta$}{Damping ratio.}
 \nomenclature[C, 30]{$\omega_n$}{Natural frequency $[\mathrm{p.u.}]$.}
 \nomenclature[C, 31]{$\tau_o$}{Weighting factor of representative day $o$.}
  \nomenclature[C, 32]{$\epsilon$}{Corrective power deviation tolerance $[\mathrm{kW}]$.}

 \nomenclature[D, 01]{$\Omega^\mathrm{C}$}{Set of CIGs, $\Omega^\mathrm{C_n}$ {\color{black}being} the set of generators connected to node $n$.}
 \nomenclature[D, 02]{$\Omega^\mathrm{C}_{d/v}$}{Set of CIGs with droop/VSM control scheme.}
 \nomenclature[D, 03]{\scalebox{0.95}{$\Omega^\mathrm{gm,opr}$}}{Feasible space of operational variables in grid-connected mode.}
 \nomenclature[D, 04]{\scalebox{0.95}{$\Omega^\mathrm{im,opr}$}}{Feasible space of operational variables in islanded mode.}
 \nomenclature[D, 05]{\scalebox{0.95}{$\Omega^\mathrm{inv}$}}{Feasible space of investment-related variables.}
 \nomenclature[D, 06]{$\Omega^\mathrm{L}$}{Set of lines connecting neighbouring nodes.}
  \nomenclature[D, 07]{$\color{black}\Omega^\mathrm{MG}$}{\color{black}Feasible space of the MG planning problem.}
 \nomenclature[D, 07]{$\Omega^\mathrm{N}$}{Set of nodes,  $\Omega^\mathrm{N_n}$ {\color{black}being} the set of nodes after and connected to node $n$.}
 \nomenclature[D, 08]{$\Omega^\mathrm{O}$}{Set of representative days.}
 \nomenclature[D, 09]{$\Omega^\mathrm{S}$}{Set of SGs, $\Omega^\mathrm{S_n}$ {\color{black}being} the set of generators connected to node $n$.}
 \nomenclature[D, 10]{$\Omega^\mathrm{T}$}{Set of hours.}

 \nomenclature[E, 01]{$\underline \bullet$/$\overline \bullet$}{Lower/upper bounds of the quantity $\bullet$.}
 \nomenclature[E, 01]{$\hat \bullet$}{Deviations of the quantity $\bullet$ in the islanded mode from its value in the grid-connected mode, {\color{black}{$\hat \bullet^+$/$\hat \bullet^-$} being upward/downward deviations}.}
 \nomenclature[E, 01]{$\mid\bullet\mid$}{Cardinality of the set $\bullet$.}

 \nomenclature[F, 01]{$d_{nto}^\mathrm{p/q}$}{Active/reactive load of node $n$ at hour $t$ in representative day $o$ $[\mathrm{kWh}/\mathrm{kVAr}]$.}
 \nomenclature[F, 01]{$d_{nto}^\mathrm{p_f/q_f}$}{Flexible part of active/reactive load $d_{nto}^\mathrm{p/q}$ $[\mathrm{kWh}/\mathrm{kVAr}]$.}
 \nomenclature[F, 01]{$p/q_{gto}$}{Active/reactive power generation of generator $g$ at hour $h$ in representative day $o$ $[\mathrm{kWh}/\mathrm{kVAr}]$.}
 \nomenclature[F, 01]{$p_{nn'to}$}{Active power flow of a line connecting nodes $(n,n')$ at hour $h$ in representative day $o$ $[\mathrm{kWh}]$.}
 \nomenclature[F, 01]{$p/q_{to}^\mathrm{b/s}$}{Active/reactive power flow bought/sold to the main grid at hour $h$ in representative day $o$ $[\mathrm{kWh}/\mathrm{kVAr}]$.}
 \nomenclature[F, 01]{$q_{nn'to}$}{Reactive power flow of a line connecting nodes $(n,n')$ at hour $h$ in representative day $o$ $[\mathrm{kVAr}]$.}
 \nomenclature[F, 01]{$v_{nto}$}{Voltage magnitude of node $n$ at hour $h$ in representative day $o$ $[V]$.}
 \nomenclature[F, 01]{$y_{nto}$}{Binary variable indicating the connection status of load of node $n$ at hour $t$ in representative day $o$ (i.e., $1/0$: connected/disconnected).}
 \nomenclature[F, 01]{$z_{g}$}{Binary variable indicating the investment status of generator $g$ (i.e., $1/0$: built/not-built).}
 \nomenclature[F, 01]{$z_{nn'}$}{Binary variable indicating the investment/reinforcement status of a line connecting nodes $(n,n')$ (i.e., $1/0$: built/not-built).}

 \nomenclature[G, 01]{$\color{black}\chi$}{\color{black}Vector of all investment and operational variables.}
 \nomenclature[G, 02]{\scalebox{0.95}{$\color{black}\chi^\mathrm{gm,opr}$}}{\color{black}Vector of ``wait-and-see'' operational variables in grid-connected mode.}
  \nomenclature[G, 03]{\scalebox{0.95}{$\chi^\mathrm{im,opr}$}}{Vector of ``wait-and-see'' operational variables in islanded mode.}
 \nomenclature[G, 04]{$\chi^\mathrm{inv}$}{Vector of ``here-and-now'' investment variables.}
 \nomenclature[G, 05]{$\eta$}{Vector of representative days (i.e., scenarios).}

 \printnomenclature

\vspace{-0.3cm}
\section{Introduction}

\lettrine[lines=2]{R}{esilient} electric networks must have the ability to ride through extreme contingencies, maintain basic service levels to critical load demands, and ensure fast recovery to normality. In other words, a resilient system should be able to modify its functionality and alter its structure in an agile manner without collapsing~\cite{Grid-Resilience-Elasticity-Is-Needed-When-Facing-Catastrophes}. The main measures to enhance the resilience of electric networks can be categorized into~\cite{Research-on-Resilience-of-Power-Systems-Under-Natural-Disasters}: (i) {\it ``hardening''}, which incorporates all activities aimed at reinforcing electric networks and enhancing component designs and constructions with the intention of preserving functionality and minimizing damage; (ii) {\it ``survivability''}, which includes innovative technologies to diversify energy supply and improve system flexibility; and (iii) {\it ``recovery''}, which incorporates all tools aimed at restoring the system to normal operating conditions. Of particular concern is the resilience of electricity distribution networks due to their interdependence with other critical infrastructure, which might culminate in a sustained negative impact on society. With this background, Microgrids (MGs) have been widely considered as a potential pathway for enhancing system resilience and ensuring both structural reinforcement and operational flexibility by allowing for the coexistence of Distributed Energy Resources (DERs) with the traditional bulk grid~\cite{Microgrids-as-a-resilience-resource-and-strategies-used-by-microgrids-for-enhancing-resilience,7091066}.

MGs are flexible distribution systems able to operate in both grid-connected and islanded mode~\cite{MicroGrids}. Their islanding capability is critical in enhancing resilience by ensuring continuity and mitigating interruptions of energy supply to consumers in the event of extreme weather conditions or significant faults in the bulk transmission grid~\cite{7489002,ZHOU2018374}. The successful island creation, especially following disastrous events, is subject to the secure transient performance of DERs, thus ensuring the survivability of the MG. However, unlike traditional bulk grids, MGs are inherently faced with a lack of rotational inertia and damping capability affecting their security in the event of significant power imbalance, and more importantly, islanding from the main grid~\cite{8450880}. Thus, it is vital to design a resilient and reliable MG able to withstand both High-Impact-Low-Frequency (HILF) and Low-Impact-High-Frequency (LIHF) uncertainties, under static and transient operational constraints. On the one hand, static islanding constraints ensure the MG's operational adequacy in supplying the forecasted electricity demand. On the other hand, transient islanding constraints ensure the MG's operational security by adhering to a dynamic response within the defined regulatory limits and, consequently, avoiding the operation of protective devices that would result in DER disconnections.

Different stochastic~\cite{yamangil2015resilient}, robust~\cite{6920097,LIN20181266,7381672}, and distributionally robust~\cite{8274133} planning tools have been previously presented in the literature for optimal investment in distribution networks and MGs aiming at enhancing system resilience to extreme contingencies. In~\cite{yamangil2015resilient}, a stochastic model has been proposed for optimal investment in distribution networks under different disastrous events characterized by a set of scenarios. In~\cite{6920097}, a robust resilience-constrained MG planning model is introduced under the uncertainty of loads and power generation of Renewable Energy Sources (RESs), with islanding from the main grid considered as another {\color{black}source of uncertainty}. In~\cite{LIN20181266}, a robust defender-attacker-defender model is presented for optimal hardening planning in resilient distribution networks, considering topology reconfiguration and islanding formation. Besides, in~\cite{7381672}, a robust model is proposed for hardening and investment planning in distribution networks based on a multi-stage and multi-zone uncertainty modeling of spatial and temporal characteristics of natural disasters. Additionally, a distributionally robust resilience-constrained investment planning model under natural disasters is introduced in~\cite{8274133}, where a moment-based ambiguity set characterizes extreme events. Even though the literature offers several operational planning models for traditional bulk grids under dynamic security constraints~\cite{5288558,867137,7232567,7115982,8269332, 6717054,UrosUC}, to the best of the authors' knowledge, previous resilience-constrained investment planning models for MGs~\cite{Microgrid-planning-considering-the-resilience-against-contingencies,6920097}, and even active distribution networks~\cite{yamangil2015resilient,6920097,LIN20181266,7381672}, have only considered static operational constraints rather than dynamic ones.  

The problem of ensuring dynamic security in power systems has been studied in~\cite{5288558}, where a transient stability-constrained Optimal Power Flow (OPF) is employed with a single-machine infinite-bus model characterizing the transient stability constraints in order to facilitate secure frequency response. Similarly, a discretized transient response is embedded in the OPF problem in~\cite{867137} to ensure a secure transient frequency response. In~\cite{7232567}, an analytical formulation is presented to limit the Rate-of-Change-of-Frequency (RoCoF) based on a single-machine Center-of-Inertia (CoI) frequency model, while~\cite{7115982} uses a simplified model of transient frequency metrics to analyze the post-fault response. Such simplified frequency response models tend to describe system dynamics inaccurately and cannot quantify the support provided by different units. 
A reduced second-order model is used in~\cite{8269332} to determine sufficient synthetic inertia and droop slopes for a collection of traditional and inverter-interfaced generators that satisfy both steady-state and dynamic frequency requirements. Moreover, in~\cite{6717054} and~\cite{UrosUC}, the unit commitment problem is solved under frequency-related constraints for traditional and low-inertia grids, where frequency-related constraints are derived based on a low-order non-linear frequency response model~\cite{65898}. 

Nevertheless, all aforementioned studies have certain drawbacks, as they are based on either simplified dynamic models~\cite{7232567,7115982,8269332}, linearized frequency-related constraints~\cite{6717054}, or make \textit{ex-ante} bound extractions on the relevant variables~\cite{UrosUC} to simplify the planning model. Furthermore, the simplifications therein represent the characteristic properties of transmission networks rather than active distribution networks and MGs. Accordingly, it is vital to present a resilience-oriented MG planning tool, including both static and transient constraints, based on a detailed dynamic model to ensure satisfactory operation given the abrupt main grid disconnection in the event of extreme contingencies.

The paper's main contributions can be summarized as follows: 
\begin{itemize}
\item A min-max-min, stochastic-robust, investment planning model, is introduced to design a resilient MG under both HILF and LIHF uncertainties. The HILF uncertainty pertains to the unscheduled islanding of the MG from the main grid while the LIHF uncertainties relate to correlated load and RES generation. For the latter, the $k$-means clustering technique is used to obtain a sufficient number of scenarios (i.e., representative days) characterizing different realizations of LIHF uncertainties.
\item Both static and transient islanding constraints (i.e., the maximum RoCoF and the frequency nadir as transient-state criteria, and the frequency deviation as a quasi steady-state criterion) are considered in the proposed model to ensure resilience under HILF and LIHF uncertainties. To the best of the authors' knowledge, there is no similar planning tool in the literature that includes both static and transient islanding constraints.
\item A tractable three-stage solution approach is presented since the proposed min-max-min, hybrid, stochastic-robust investment problem with a non-linear frequency response model cannot be solved directly.
\end{itemize}

The rest of the paper is organized as follows. In Section \ref{sec: assumptions}, the investment planning model is described in a compact form together with the main modeling preliminaries. Section~\ref{sec: formulation} presents {\color{black}the} detailed problem formulation under static and transient islanding constraints, whereas Section~\ref{sec: case studies} discusses the application of the proposed investment planning model on {\color{black}the} CIGRE 18-node distribution network in order to design a resilient MG under different operating conditions. Finally, Section~\ref{sec: conclusion} concludes the paper.  

\vspace{-0.25cm}
\section{Modeling Preliminaries and Problem Description}\label{sec: assumptions}
The modeling preliminaries in the proposed investment planning model are: 
\begin{itemize}
    \item Without loss of generality, a single-year planning horizon is considered rather than a multi-year one to reach a compromise between accuracy and tractability of the proposed model.  
    \item To capture interday/intraday variation/ramping of uncertain loads and power generation of RESs, a sufficient number of representative days (i.e., {\it scenarios}) is considered, obtained by the $k$-means clustering technique~\cite{dehghan2019robust}. 
    \item A linearized version of the DistFlow model is used for the power flow formulation to obtain a linear optimization problem~\cite{Network-reconfiguration-in-distribution-systems-for-loss-reduction-and-load-balancing, Robust-Optimization-Based-Optimal-DG-Placement-in-Microgrids}. Additionally, the quadratic line flow expressions are linearized using a piecewise linear approximation~\cite{A-Linearized-OPF-Model-With-Reactive-Power-and-Voltage-Magnitude}. Finally, a constant marginal cost is utilized to eliminate the non-linearity of quadratic cost functions~\cite{dehghan2019robust}.
    \item It is assumed that an unscheduled islanding event might happen at \textit{each} hour of the representative days considered.
\end{itemize}

\vspace{-0.35cm}
\subsection{Compact Formulation under Static Constraints}\label{sec: compact formulation}
The proposed min-max-min investment planning model under static operational constraints in grid-connected and islanded mode can be presented in compact form {\color{black}as:}
\begin{align}\label{eqn: compact}
\begin{split}
        \min_{\color{black}\chi \in \Omega^\mathrm{MG}}\;\,\Theta^{{\mathrm{inv}}}(\chi^{\mathrm{inv}})&+\Theta^{\mathrm{gm,opr}}(\chi^{\mathrm{inv}},\chi^{\mathrm{gm,opr}})\\
        &+||\breve{\Theta}^{\mathrm{im,opr}}(\chi^{\mathrm{inv}},\chi^{\mathrm{gm,opr}},\chi^{\mathrm{im,opr}})||_{\infty},
\end{split}
\end{align}
where {\color{black} $\Omega^{\mathrm{MG}} = \{ \chi = [\chi^\mathrm{inv}, \chi^\mathrm{gm,opr},\chi^\mathrm{im,opr}] \,|\, \chi^{\mathrm{inv}} \in \Omega^{\mathrm{inv}}\,;\,\chi^{\mathrm{gm,opr}}\in \Omega^{\mathrm{gm,opr}}\,;\,\chi^{\mathrm{im,opr}}\in\Omega^{\mathrm{im,opr}}$\}}, $\breve{\Theta}^{\mathrm{im,opr}}=[\min\,\Theta_{11}^{\mathrm{im,opr}},...,\min\,\Theta_{\mathcal TO}^{\mathrm{im,opr}}]$, $\color{black}{\mathcal T}=|\Omega^{\mathrm{T}}|$, and $\color{black}{\mathcal O} =| \Omega^{\mathrm{O}}|$. Also, $||\breve{\Theta}^{\mathrm{im,opr}}||_{\infty} = \max (\underset{\forall t,\forall o}\min \,\Theta_{to}^{\mathrm{im,opr}})$. Hence, the objective function \eqref{eqn: compact} minimizes the total investment costs ($\Theta^{{\mathrm{inv}}}$), the {\it ``expected''} total operation costs in grid-connected mode for all hours of all representative days ($\Theta^{{\mathrm{gm,opr}}}$), and the {\it ``worst-case''} total penalty costs of disconnecting loads from MG in islanded mode for all hours in all representative {\color{black}days} ($\Theta^{{\mathrm{im,opr}}}$). 
 
The min-max-min objective function \eqref{eqn: compact} can be rewritten as a single minimization problem by using the auxiliary variable $\gamma$:
\begin{subequations}
\label{eqn:com}
\begin{alignat}{3}
    & \min_{\color{black}\chi \in \Omega^{\mathrm{MG}}} \; && \Theta^{{\mathrm{inv}}}(\chi^{\mathrm{inv}})+\Theta^{\mathrm{gm,opr}}(\chi^{\mathrm{inv}},\chi^{\mathrm{gm,opr}})+\gamma \label{eqn:com1}\\
    & \;\mathrm{s.t.} \;&& \gamma\geq\Theta_{to}^{\mathrm{im,opr}}(\chi^{\mathrm{inv}},\chi^{\mathrm{gm,opr}},\chi^{\mathrm{im,opr}}), \; \forall t\in\Omega^\mathrm{T},o\in\Omega^\mathrm{O}, \label{eqn:com 2}
\end{alignat}
\end{subequations}
The optimization problem \eqref{eqn:com} is a Mixed-Integer Linear Programming (MILP) problem, and as such can be solved by available software packages to obtain optimal investment and operation decisions in grid-connected and islanded mode. However, the operation decisions may violate transient islanding constraints. To remedy such limitation and ensure MG resilience before and after an islanding event, a three-stage methodology is employed, incorporating a non-linear model for evaluation of the transient frequency response of a MG after islanding. In the sequel, the MG frequency dynamics, the metrics to evaluate the transient frequency response of a MG in islanded mode, as well as the proposed three-stage solution approach are presented.
\vspace{-0.5cm}
\subsection{Microgrid Frequency Dynamics}\label{sec: MG dynamics}

The employed dynamic model is based on the uniform representation of frequency transients in a low-inertia system previously introduced in~\cite{UrosLQR,UrosUC}, comprising both traditional Synchronous Generators (SGs, indexed by $i \in \Omega^\mathrm{S}$) and Converter-Interfaced Generators (CIGs, indexed by $c \in \Omega^\mathrm{C}$). The generator dynamics are described by the CoI swing equation with aggregate inertia $M_g$ and damping $D_g$. The low-order model proposed in~\cite{Anderson1990} is used for modeling the governor droop and turbine dynamics.
The impact of {\it grid-supporting} CIGs providing frequency support via droop $(d\in \Omega^\mathrm{C}_d\subseteq\Omega^\mathrm{C})$ and virtual synchronous machine (VSM) $(v \in \Omega^\mathrm{C}_v\subseteq\Omega^\mathrm{C})$ control is also included, as these are the two most common control approaches in the literature~\cite{Rocabert2012,UrosGM}. 
Hence, the transfer function $G(s)$ between the active power change $\Delta P_e(s)$, {\color{black}with positive values corresponding to a net load decrease,} and the CoI frequency deviation $\color{black}\Delta f(s)$ can be derived as:   
\begin{align}
    G(s) &= \dfrac{\Delta f(s)}{\Delta P_e(s)} = \Bigg(\underbrace{(sM_g+D_g )+\sum\limits_{i\in\Omega^\mathrm{S}} \dfrac{K_{i} (1+sF_{i} T_{i} )}{R_{i}(1+sT_{i})}}_{\text{traditional SGs}} \nonumber \\
    & + \underbrace{\sum\limits_{d \in \Omega^\mathrm{C}_d} \dfrac{K_{d}}{R_{d} (1+sT_{d} )}}_{\text{droop-based CIGs}} + \underbrace{\sum\limits_{v\in\Omega^\mathrm{C}_v } \dfrac{sM_{v}+D_{v}} {1+sT_{v}} }_{\text{VSM-based CIGs}}\Bigg)^{-1}. \label{eq:G1}
\end{align}
Assuming that the time constants $(T_i\approx T)$ of all SGs are several orders of magnitude higher than the ones of converters~\cite{UrosStab}, one can approximate $T \gg T_{d,v}\approx0$, which transforms \eqref{eq:G1} into:
\begin{equation}
    G(s) = \frac{1}{MT}\frac{1+sT}{s^2+2\zeta\omega_n s + \omega_n^2}, \label{eq:G2}
\end{equation}
where $\omega_n = \sqrt{\frac{D+R_g}{MT}}$ and $\zeta = \frac{M+T(D+F_g)}{2\sqrt{MT(D+R_g)}}$. 
More details on the proposed second-order frequency model in \eqref{eq:G2} and mathematical formulation can be found in~\cite{UrosLQR}.

\vspace{-0.25cm}
\subsection{Dynamic Metrics for Microgrid Islanding} \label{sec:metrics}
Following a disturbance, the dynamic frequency response is characterized by the instantaneous RoCoF ($\dot{f}_\mathrm{max}$) and frequency nadir ($\Delta f_\mathrm{max}$), whereas the steady-state response is governed by the constant frequency deviation from a pre-disturbance equilibrium ($\Delta f_\mathrm{ss}$). By assuming a stepwise disturbance in the active power $\Delta P_e(s) = -\Delta P/s$, {\color{black} where $\Delta P$ is the net power change,} the time-domain expression for frequency metrics of interest can be derived as follows:
\begin{subequations} \label{eq:puConst}
\begin{align}
    \dot{f}_\mathrm{max} &= \dot{f}(t_0^+) = -\frac{\Delta P}{M}, \label{eq:rocof}\\
    \Delta f_\mathrm{max} &= - \frac{\Delta P}{D+R_g} \left( 1 + \sqrt{\dfrac{T(R_g-F_g)}{M}} e^{-\zeta\omega_n t_m} \right), \label{eq:nadir}\\
    \Delta f_\mathrm{ss} &= -\frac{\Delta P}{D + R_g}, \label{eq:qss}
\end{align}
\end{subequations}
with the introduction of new variable $\omega_d = \omega_n\sqrt{1-\zeta^2}$ and $t_m =(\sfrac{1}{\omega_d})\tan^{-1}\left(\sfrac{\omega_d}{\left(\omega_n\zeta - T^{-1}\right)}\right)$ denoting the time instance of frequency nadir.

It can be clearly seen that the aggregate system parameters such as $M$, $D$, $R_g$ and $F_g$ have a direct impact on frequency performance. In particular, RoCoF and steady-state deviation are explicitly affected by $M$ and $(D,R_g)$, respectively, while frequency nadir has a non-linear dependency on all four system factors. With the increasing penetration of CIGs and subsequent decommissioning of conventional SGs, these parameters are drastically reduced and can compromise the overall frequency performance. To prevent the accidental activation of load-shedding, under/over frequency and RoCoF protection relays, the proposed three-stage solution algorithm, described in {\color{black} the following}, imposes limits on the aforementioned frequency metrics to account for low levels of inertia and damping and their impact on the frequency response after a MG islanding.

\vspace{-0.35cm}
\subsection{Three-Stage Solution Algorithm}\label{sec: solution algorithm}

\begin{figure}[t!] 
	\centering
    \scalebox{0.98}{\includegraphics[width=\linewidth]{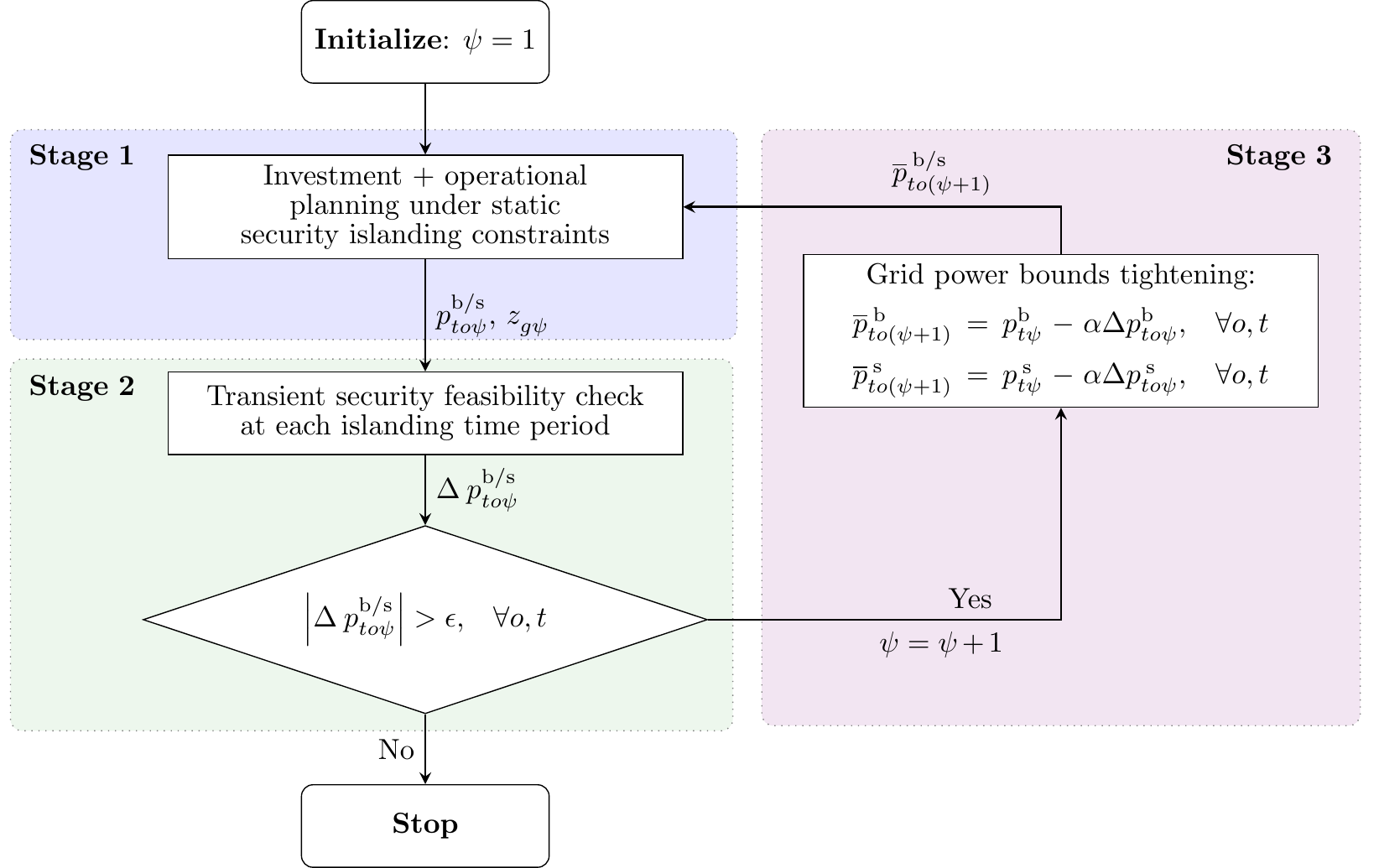}}
	\caption{Proposed three-stage MG planning algorithm.}
	\label{fig:methodology}
	\vspace{-0.45cm}
\end{figure}

The proposed three-stage approach proposed for solving the problem \eqref{eqn:com} with the inclusion of static and transient security constraints is illustrated in~Fig.~\ref{fig:methodology} and can be summarized as follows: 
\vspace{+0.2em}
\subsubsection{Solving the Static Investment Planning Problem} 
At each iteration $\psi$, the investment planning model \eqref{eqn:com} is solved under static security constraints in both grid-connected and islanded mode. A detailed formulation is provided in Section~\ref{sec:extended}.
\vspace{+0.2em}
\subsubsection{Evaluating {\color{black}{Transient Frequency Security}}}
The unscheduled loss of power exchange with the main grid may result in large frequency transients within the MG. Following the discussion from Section~\ref{sec:metrics}, the transient frequency response is characterized by the magnitude of the abrupt active power loss ($\Delta P$) and the aggregate control parameters of all MG generators ($M$, $D$, $R_g$, and $F_g$). Therefore, the magnitude of the potential disturbance at each hour of every representative day during the planning horizon is set equal to the power exchange with the main grid scheduled at the time of the disconnection (i.e., $\Delta P = p^\mathrm{b/s}_{to\psi}$). At the same time, the control parameters depend on the built/not-built status of generators in the MG at each iteration ($z_{g\psi}$). Consequently, at each iteration $\psi$, this stage of the algorithm uses the variables $p^\mathrm{b/s}_{to\psi}$ and $z_{g\psi}$ obtained from the first stage to evaluate the solution feasibility under transient security constraints in \eqref{eq:puConst}, described in detail in Section~\ref{sec:dy}. 

The solution of the second stage provides the minimum amount of {\color{black}corrective deviation ($\Delta p^\mathrm{b/s}_{to\psi}$) from the scheduled power exchange with the main grid} ($p^\mathrm{b/s}_{to\psi}$) to meet the transient security criteria. If {\color{black}this value is zero or less than a small tolerance,} the optimal investment and operational solution obtained from the first-stage problem ensures frequency security in islanded mode, and the algorithm is terminated. Otherwise, the algorithm proceeds to the third stage.

\vspace{+0.2em}
\subsubsection{Tightening Power Exchange with the Main Grid}
The third stage of the algorithm employs the non-zero solution obtained from the previous stage to tighten the permissible limits imposed on the power exchange with the main grid at each hour of every representative day throughout the planning horizon. Once the respective bounds have been altered, the algorithm proceeds to the next iteration ($\psi+1$). The modified limits may lead either to a change in the power exchange with the main grid (through operational decisions), a change in the investment decisions, or both.

\vspace{-0.35cm}
\section{Extended Formulation under Static and Transient Constraints}\label{sec: formulation}
In this section, the extended formulation of the proposed planning tool under static and transient constraints is presented. The iteration index $\psi$ is omitted for better legibility and brevity.

\vspace{-0.25cm}
\subsection{\color{black}Extended Formulation under Static Constraints}\label{sec:extended}
{\color{black}The optimization problem in the first stage corresponds to the compact formulation \eqref{eqn:com}, including investment limitations and static operational constraints in grid-connected and islanded mode. The extended terms in the objective function and the constraints are outlined in the following.}  
\vspace{+0.3em}
\subsubsection{\color{black}Investment}\label{sec:inv}
The {\color{black}term} $\Theta^{\mathrm{inv}}(\chi^{\mathrm{inv}})$ {\color{black}in the objective function} of the proposed planning problem is given by: 
\begin{equation}\label{1}
   \Theta^{\mathrm{inv}} = \sum_{g\in\left\{ \Omega^{\mathrm{S}},\Omega^{\mathrm{C}}\right\} }\left(ic_{g}\cdot z_{g}\right)+\sum_{(n,n')\in\Omega^{\mathrm{L}}}\left(ic_{nn'}\cdot z_{nn'}\right),
\end{equation}
and includes the total investment/reinforcement costs of generators/lines throughout the planning horizon. The  optimization variables $\chi^{\mathrm{inv}}=\{z_g,z_{nn'}\},\forall g\in\left\{ \Omega^{\mathrm{S}},\Omega^{\mathrm{C}}\right\} \land \forall (n,n')\in\Omega^{\mathrm{L}}$ are {\it here-and-now} decisions. 

\vspace{+0.3em}
\subsubsection{\color{black}Grid-Connected Operation}\label{sec:gc}
The {function $\Theta^{\mathrm{gm,opr}}$ capturing the operational cost in grid-connected mode is given by:} 
\begin{subequations}
    \begin{align}
        \begin{split}
            \Theta^{\mathrm{gm,opr}} = &\sum_{o\in\Omega^{\mathrm{O}}}\sum_{t\in\Omega^{\mathrm{T}}}\left(\tau_{o}\cdot\left(e_{to}^{\mathrm{b}}\cdot p_{to}^{\mathrm{b}}-e{}_{to}^{\mathrm{s}}\cdot p_{to}^{\mathrm{s}}\right)\right)\\
            +&\sum_{o\in\Omega^{\mathrm{O}}}\sum_{t\in\Omega^{\mathrm{T}}}\sum_{g\in\left\{ \Omega^{\mathrm{S}},\Omega^{\mathrm{C}}\right\} }\left(\tau_{o}\cdot mc_{g}\cdot p_{gto}\right)\\
            +&\sum_{o\in\Omega^{\mathrm{O}}}\sum_{t\in\Omega^{\mathrm{T}}}\sum_{n\in \Omega^{N} }\left(\tau_{o}\cdot fc_{n}\cdot d_{nto}^{\mathrm{p}_f}\right)\label{2a}
        \end{split}
    \end{align}
    {\color{black}The constraints that need to be taken into account to reflect operational limitations in grid-connected mode are:} 
    \begin{align}
        \begin{split}
            &p_{n^{''}nto}+ {\color{black}p_{to|n=1}^\mathrm{b}-p_{to|n=1}^\mathrm{s}} +\sum_{g\in\{\Omega^\mathrm{S_{n}},\Omega^\mathrm{C_{n}}\}}p_{gto}\\ 
            &\quad =\sum_{n^{\prime}\in\Omega^\mathrm{N_{n}}}p_{nn^{\prime}to}+d_{nto}^\mathrm{p}, \quad \forall n\in\Omega^\mathrm{N},t\in\Omega^\mathrm{T},o\in\Omega^\mathrm{O}, \label{2b}
        \end{split}\\
        \begin{split}
            &q_{n^{''}nto}+{\color{black}q_{to|n=1}^\mathrm{b}-q_{to|n=1}^\mathrm{s}}+\sum_{g\in\{\Omega^\mathrm{S_{n}},\Omega^\mathrm{C_{n}}\}}q_{gto}\\ 
            &\quad =\sum_{n^{\prime}\in\Omega^\mathrm{N_{n}}}q_{nn^{\prime}to}+ d_{nto}^\mathrm{q},\quad  \forall n\in\Omega^\mathrm{N},t\in\Omega^\mathrm{T},o\in\Omega^\mathrm{O}, \label{2c}
        \end{split}\\
        \begin{split}
        &v_{n^{\prime\prime}to}-v_{nto}=
        \left(r_{n^{\prime\prime}n}\cdot p_{n^{\prime\prime}nto}+x_{n^{\prime\prime}n}\cdot q_{n^{\prime\prime}nto}\right), \\ 
        & \qquad \qquad \qquad \qquad \qquad \qquad \,\, \forall n\in\Omega^\mathrm{N},t\in\Omega^\mathrm{T},o\in \Omega^\mathrm{O}, \label{2d}
        \end{split}\\
        \begin{split}
        &0 \leq p_{to}^\mathrm{b} \leq {\overline p}_{to}^\mathrm{b}, \; 0 \leq p_{to}^\mathrm{s} \leq \overline {p}_{to}^\mathrm{s}, \quad \quad \quad \,\,\, \forall t\in\Omega^\mathrm{T},o \in \Omega^\mathrm{O}, \label{2f}
        \end{split}\\
        \begin{split}
        &{\color{black} 0 \leq q_{to}^\mathrm{b} \leq {\overline q}_{to}^\mathrm{b}, \; 0 \leq q_{to}^\mathrm{s} \leq \overline {q}_{to}^\mathrm{s}, \quad \quad \quad \,\,\,\, \forall t\in\Omega^\mathrm{T},o \in \Omega^\mathrm{O},} \label{2ff}
        \end{split}\\
        &d_{nto}^{\mathrm{p}}=d_{nto}^{\mathrm{p_{c}}}+d_{nto}^\mathrm{p_{f}}, \quad \quad \quad \quad \, \forall n\in\Omega^\mathrm{N},t\in\Omega^\mathrm{T},o\in\Omega^\mathrm{O}, \label{2g}\\
       & d_{nto}^\mathrm{q}=d_{nto}^\mathrm{q_{c}}+d_{nto}^\mathrm{q_{f}},\quad \quad \quad \quad \, \forall n\in\Omega^\mathrm{N},t\in\Omega^\mathrm{T},o\in\Omega^\mathrm{O}, \label{2h}\\
        &\underline{d}_{nto}^\mathrm{p_{f}}\leq d_{nto}^\mathrm{p_f}\leq\overline{d}_{nto}^\mathrm{p_{f}}, \quad \quad \quad \quad \forall n\in\Omega^\mathrm{N},t\in\Omega^\mathrm{T},o\in\Omega^\mathrm{O}, \label{2i}\\
        &\underline{d}_{nto}^\mathrm{q_{f}}\leq d_{nto}^\mathrm{q_f}\leq\overline{d}_{nto}^\mathrm{q_{f}}, \quad \quad \quad \quad \forall n\in\Omega^\mathrm{N},t\in\Omega^\mathrm{T},o\in\Omega^\mathrm{O}, \label{2j}\\
        &\sum_{t\in\Omega^{T}}{\color{black}d_{nto}^\mathrm{p}}=e_{no}, \quad \quad \quad \quad \quad \quad \quad \quad \,\,\, \forall n\in\Omega^\mathrm{N},o\in\Omega^\mathrm{O}, \label{2k}\\
        \begin{split}
        &0\leq p_{gto}\le {\overline p}_{gto}\cdot z_{g}, \quad \, \forall g\in \{\Omega^\mathrm{C},\Omega^\mathrm{S}\},t\in\Omega^\mathrm{T},o \in \Omega^\mathrm{O}, \label{2l}
        \end{split}\\
        \begin{split}
        &{\underline q}_{gto}\cdot z_{g}\leq q_{gto}\le {\overline q}_{gto}\cdot z_{g}, \quad \forall  g\in\{\Omega^\mathrm{C},\Omega^\mathrm{S}\},t\in\Omega^\mathrm{T},o \in \Omega^\mathrm{O}, \label{2m}
        \end{split}\\
        &-r_{g}^\mathrm{d}\leq p_{gto}-p_{g(t-1)o}\leq r_{g}^\mathrm{u},\forall g\in\Omega^\mathrm{S},t\in\Omega^\mathrm{T},o\in\Omega^\mathrm{O}, \label{2n}\\
        &\sum_{t\in\Omega^{T}}p_{gto}\leq {\color{black}c_{go}\cdot {\overline p}_{gto}}, \quad \quad \quad\,\forall g\in\Omega^\mathrm{S},t\in\Omega^\mathrm{T},o\in\Omega^\mathrm{O}, \label{2o}\\
        \begin{split}
        &p_{nn'to}^{2}+q_{nn'to}^{2}\leq \overline{s}_{nn'}^{2}\cdot\left(z_{nn'}^{0}+z_{nn'}\right),\\ 
        &\qquad \qquad \qquad \qquad \qquad \;\, \forall (n,n')\in\Omega^\mathrm{L},t\in\Omega^\mathrm{T},o\in\Omega^\mathrm{O}, \label{2p}
        \end{split}\\
        \begin{split}
        &\underline v\le v_{nto}\le \overline v, \; v_{to|n=1}=1, \;\,\,\,  \forall n\in\Omega^\mathrm{N},t\in\Omega^\mathrm{T},o\in\Omega^\mathrm{O},\label{2q}
        \end{split}
\end{align}
\end{subequations}\label{2}
Here, the vector of {\it wait-and-see} decision variables is given by $\chi^{\mathrm{gm,opr}}=\{d^{\mathrm{p/q}}_{nto}, d^{\mathrm{p_f/q_f}}_{nto}, p_{gto}, p_{nn'to}, p^{\mathrm{b/s}}_{to}, q_{gto}, q_{nn'to}, v_{nto}\}$.

The objective function \eqref{2a} minimizes the total operation costs, including the total costs of power exchange with the main grid, the total operation costs of generators, and the total penalty costs of shifting loads away from the periods preferred by consumers. Constraints \eqref{2b}-\eqref{2d} describe the power flows based on the linearized version of the DistFlow model~\cite{Network-reconfiguration-in-distribution-systems-for-loss-reduction-and-load-balancing, Robust-Optimization-Based-Optimal-DG-Placement-in-Microgrids}, and \eqref{2f}{\color{black}-\eqref{2ff}} ensure the non-negativity and impose the upper limits on the power exchange with the main grid. Furthermore, \eqref{2g}-\eqref{2j} reflect the power balance of constant and flexible loads as well as the limitations of flexible loads at each node and at every hour of each representative day, whereas \eqref{2k} ensures that the daily energy consumption of flexible loads is maintained for each representative day. Constraints \eqref{2l}-\eqref{2n} denote capacity and ramp-rate limits of generators at each hour of every representative day. Moreover, the reactive power limits of CIGs are based on the maximum generated active power, i.e., $\overline{q}^{}_{gto} = \tan{\phi}\cdot \overline{p}^{}_{gto}$, where $\cos{\phi}$ is the maximum power factor of a unit defined by the grid code. Constraint \eqref{2o} defines the daily capacity factor of SGs in each representative day of the planning horizon{\color{black}\cite{turvey1977electricity}}, and \eqref{2p} imposes the thermal loading limits of each line. The latter quadratic constraint is linearized by means of a convex polygon, defined by inner approximations of the thermal loading circle~\cite{A-Linearized-OPF-Model-With-Reactive-Power-and-Voltage-Magnitude}. 
Finally, \eqref{2q} 
limits the nodal voltage magnitudes throughout the planning horizon. 
\vspace{+0.3em}
\subsubsection{\color{black}Islanded Operation}\label{sec:im}
It is assumed that at every hour of each representative day, the MG should be able to withstand an unscheduled islanding event. The operation planning problem of a MG in islanded mode is aimed at ensuring survivability and self-sufficiency, where priority is given to critical loads. It is worthwhile to note that, in this paper, the self-sufficiency is ensured for one period (i.e., one hour) after disconnection from the main grid. However, the islanded operation period can be straightforwardly extended to multiple periods based on the required resilience. Hereafter, the superscript $\textrm {``im''}$ denotes operational variables in islanded mode. The {function $\Theta_{to}^{\mathrm{im,opr}}$ capturing the operational cost in islanded mode is given by:}
\begin{subequations}\label{3}
    \begin{align}\label{3a}
    \begin{split}
         \Theta_{to}^{\mathrm{im,opr}} = \sum_{n\in\Omega^\mathrm{N}}\left(pc_{n}\cdot\left(\left(1-y_{nto}\right)\cdot d_{nto}^\mathrm{p_{c}}+ {\color{black}\hat d_{nto}^\mathrm{p_{f}^-}}\right)\right)
    \end{split}
    \end{align}
    {\color{black}The constraints that need to be taken into account to reflect operational limitations in islanded mode are:} 
    \begin{align}
        \begin{split}            
            &p_{n^{''}nto}^\mathrm{im}+\sum_{g\in\{\Omega^\mathrm{S_{n}},\Omega^\mathrm{C_{n}}\}}p_{gto}^\mathrm{im}=\sum_{n^{\prime}\in\Omega^\mathrm{N_{n}}}p_{nn^{\prime}to}^\mathrm{im}+\\
           &\left(y_{nto}\cdot d_{nto}^\mathrm{p_{c}}+d_{nto}^\mathrm{im,p_{f}}\right), \quad \quad \, \forall  n\in\Omega^\mathrm{N},t\in\Omega^\mathrm{T},o\in\Omega^\mathrm{O},\label{3b}
        \end{split}\\
        \begin{split}            
            &q_{n^{''}nto}^\mathrm{im}+\sum_{g\in\{\Omega^\mathrm{S_{n}},\Omega^\mathrm{C_{n}}\}}q_{gto}^\mathrm{im}=\sum_{n^{\prime}\in\Omega^\mathrm{N_{n}}}q_{nn^{\prime}to}^\mathrm{im}+\\
           &\left(y_{nto}\cdot d_{nto}^\mathrm{q_{c}}+d_{nto}^\mathrm{im,q_{f}}\right), \quad \quad \, \forall  n\in\Omega^\mathrm{N},t\in\Omega^\mathrm{T},o\in\Omega^\mathrm{O},\label{3c}
        \end{split}\\
        \begin{split}
            &v_{n^{\prime\prime}to}^\mathrm{im}-v_{nto}^\mathrm{im}=\left(r_{n^{\prime\prime}n}\cdot p_{n^{\prime\prime}nto}^\mathrm{im}+x_{n^{\prime\prime}n}\cdot q_{n^{\prime\prime}nto}^\mathrm{im}\right),\,\\ 
            &\,\qquad \qquad \qquad \qquad \qquad \qquad \forall n\in\Omega^\mathrm{N},t\in\Omega^\mathrm{T},o\in \Omega^\mathrm{O},\label{3d}
        \end{split}\\
            &d_{nto}^\mathrm{im,p_f}=d_{nto}^\mathrm{p_{f}}+{\color{black}\hat d_{nto}^\mathrm{p_{f}^+}-\hat d_{nto}^\mathrm{p_{f}^-}}, \,\, \forall  n\in\Omega^\mathrm{N},t\in\Omega^\mathrm{T},o\in\Omega^\mathrm{O},\label{3e}\\
            &{\color{black}0 \leq \hat d_{nto}^\mathrm{p_{f}^+},\hat d_{nto}^\mathrm{p_{f}^-}\geq 0}, \quad \quad \quad \quad \, \forall  n\in\Omega^\mathrm{N},t\in\Omega^\mathrm{T},o\in\Omega^\mathrm{O},\label{3ee}\\
            &d_{nto}^\mathrm{{im},q_{f}}=d_{nto}^\mathrm{q_{f}}+ \hat d_{nto}^\mathrm{q_{f}}, \quad \quad \quad \,\, \forall  n\in\Omega^\mathrm{N},t\in\Omega^\mathrm{T},o\in\Omega^\mathrm{O},\label{3f}\\
            &\underline{d}_{nto}^\mathrm{p_{f}}\leq d_{nto}^\mathrm{{im},p_f}\leq\overline{d}_{nto}^\mathrm{p_{f}}, \quad \quad \quad \,\, \forall  n\in\Omega^\mathrm{N},t\in\Omega^\mathrm{T},o\in\Omega^\mathrm{O},\label{3g}\\
            &\underline{d}_{nto}^\mathrm{q_{f}}\leq d_{nto}^\mathrm{{im},q_f}\leq\overline{d}_{nto}^\mathrm{q_{f}}, \quad \quad \quad \,\, \forall  n\in\Omega^\mathrm{N},t\in\Omega^\mathrm{T},o\in\Omega^\mathrm{O}, \label{3h}\\
            &{\color{black}d_{nto}^\mathrm{{im},p}}\leq e_{no}-\sum_{t'=1}^{t-1}{\color{black}d_{nt'o}^\mathrm{{im},p}}, \quad \quad \, \forall  n\in\Omega^\mathrm{N},t\in\Omega^\mathrm{T},o\in\Omega^\mathrm{O},\label{3i}\\
        \begin{split}
            &0\leq p_{gto}^\mathrm{im}\le {\overline p}_{gto}\cdot z_{g}, \quad \,\, \forall  g\in \{\Omega^\mathrm{C},\Omega^\mathrm{S}\},t\in\Omega^\mathrm{T},o \in \Omega^\mathrm{O}, \label{3j}
        \end{split}\\
        \begin{split}
            &{\underline q}_{gto}\cdot z_{g}\leq q_{gto}^\mathrm{im}\le {\overline q}_{gto}\cdot z_{g}, \quad \forall  g\in\{\Omega^\mathrm{C},\Omega^\mathrm{S}\},t\in\Omega^\mathrm{T},o \in \Omega^\mathrm{O},\label{3k}
        \end{split}\\
            &-r_{g}^\mathrm{d}\leq p_{gto}^\mathrm{im} - p_{gto}\leq r_{g}^\mathrm{u}, \quad \,\,\,\, \forall  g\in\Omega^\mathrm{S},t\in\Omega^\mathrm{T},o\in\Omega^\mathrm{O},\label{3n}\\
            &p_{gto}^\mathrm{im}\leq {\color{black}c_{go}\cdot {\overline p}_{gto}}-\sum_{t'=1}^{\color{black}t-1}p_{gto}^\mathrm{im}, \,\,\,\, \forall  g\in\Omega^\mathrm{S},t\in\Omega^\mathrm{T},o\in\Omega^\mathrm{O}, \label{3o}\\
        \begin{split}
            &{p_{nn'to}^{\mathrm{im}}}^2+{q_{nn'to}^\mathrm{im}}^2\leq s_{nn'to}^{2}\cdot\left(z_{nn'}^{0}+z_{nn'}\right),\,\\
            &\qquad \qquad \qquad \qquad \qquad \,\,\, \forall (n,n') \in\Omega^\mathrm{L},t\in\Omega^\mathrm{T},o\in\Omega^\mathrm{O}, \label{3p}
        \end{split}\\
        \begin{split}
            \underline v\le v_{nto}^\mathrm{im}\le \overline v, \; v_{to|n=1}^\mathrm{im}=1,  \quad \forall  n\in\Omega^\mathrm{N},t\in\Omega^\mathrm{T},o\in\Omega^\mathrm{O},\label{3q}
        \end{split}
\end{align}
\end{subequations}
where, similarly to the previous operation planning problem in grid-connected mode, all operation variables $\chi^{\mathrm{im,opr}}_{to}=\{{\hat d}^{\mathrm{p_f/q_f}}_{nto},d^{\mathrm{im,p_f/im,q_f}}_{nto}, {\hat p}_{gto}, p_{gto}^\mathrm{im}, p_{nn'to}^\mathrm{im}, {\hat q}_{gto}, q_{gto}^\mathrm{im}, q_{nn'to}^\mathrm{im}, v_{nto}^\mathrm{im}\}$ are {\it wait-and-see} decisions.

The objective function \eqref{3a} minimizes the total unserved load and ensures an adequate supply of at least the critical MG loads. It should be noted that $pc_{n}$ describes the priority level of the load at a specific node, with higher values suggesting more critical loads, and the amount of unserved flexible load is denoted by $\hat d_{nto}^\mathrm{p_{f}}$. Constraints \eqref{3b}-\eqref{3d} enforce the post-islanding power flow balance, whereas the deviations between the amount of flexible load served in grid-connected and islanded mode are given by \eqref{3e}-\eqref{3f} and used to determine the fraction of served and unserved flexible loads in islanded mode. Moreover, \eqref{3g}-\eqref{3h} enforce the limitations of flexible loads in islanded mode, and \eqref{3i} restricts the supply of flexible loads in terms of respective demand already served before the current time instance affected by a disconnection from the main grid. Constraints \eqref{3j}-\eqref{3k} denote capacity limits of generators, \eqref{3n}-\eqref{3o} indicate that re-scheduling actions of SGs in islanded mode are subject to their ramp rate and daily capacity factor limitations as well as their scheduling actions before the current time step. Furthermore, similar to the formulation in grid-connected mode, \eqref{3p} defines the thermal loading limit of each line and \eqref{3q} 
 limit the nodal voltage magnitudes. Note that in the grid-connected mode the voltage at the Point-of-Common Coupling (PCC) is maintained by the stiff grid, while in the islanded mode it is controlled by the DERs. 

{\color{black}The final optimization problem is a MILP problem} in the first stage of the algorithm, and its solution is subsequently used in the feasibility check in the second stage. 

\vspace{-0.25cm}
\subsection{{\color{black}Transient Security} Feasibility Checking Problem}\label{sec:dy}
The feasibility of the planning solution under transient security constraints is necessary to guarantee the secure islanding of a MG. According to the metrics described in \eqref{eq:puConst} and the discussions in Section~\ref{sec: solution algorithm}, the transient frequency response in the event of islanding depends on the amount of power exchange with the main grid at the time of {\color{black}the} event (i.e., $\Delta P = p_{to\psi}^\mathrm{b/s}$) and the control parameters of the online generators in the MG (i.e., $M(z_{g\psi})$, $D(z_{g\psi})$, $F_g(z_{g\psi})$ and $R_g(z_{g\psi})$). Note however that, with respect to decision variables, \eqref{eq:rocof} and \eqref{eq:qss} are linear while \eqref{eq:nadir} is highly non-linear. Given the optimal values of decision variables obtained from the first stage ($p^\mathrm{b/s}_{to\psi}$ and $z^{}_{g\psi}$), the non-linear term in \eqref{eq:nadir} can be defined as a constant at each iteration. Consequently, at each iteration $\psi$, the feasibility check can be formulated as a linear programming problem of the form
\begin{subequations}\label{dyn_model}
\begin{align}
 \Theta^\mathrm{dyn}_{t}  =\min_{\Delta p^\mathrm{b/s}_{to\psi}} \; \left|\Delta p^\mathrm{b/s}_{to\psi}\right|  \label{eqn:1d}
\end{align}\vspace{-0.05cm}
{\color{black}{The constraints that need to be taken into account to ensure transient security feasibility are: }}
\begin{align}
\begin{split}\label{eqn:3d}
\left | \dfrac{p^\mathrm{b/s}_{to\psi} +  \Delta p^\mathrm{b/s}_{to\psi}}{M}\right| \leq \dot{f}_\mathrm{lim}, 
\end{split}\\
\begin{split}\label{eqn:2d}
\left |\dfrac{p^\mathrm{b/s}_{to\psi} + \Delta p^\mathrm{b/s}_{to\psi}}{D + R_{g}} \cdot \left(1 + \sqrt{\dfrac{T(R_{g} - F_{g})}{M}}e^{-\zeta \omega_{n}t_{m}}\right)\right| \leq \Delta f_\mathrm{lim},
\end{split}\\
\begin{split}\label{eqn:4d}
\left | \dfrac{p^\mathrm{b/s}_{to\psi} +  \Delta p^\mathrm{b/s}_{to\psi}}{D + R_{g}}\right| \leq \Delta f_\mathrm{ss,lim}. 
\end{split}
\end{align}
\end{subequations}

The feasibility checking problem \eqref{dyn_model} is solved independently {\color{black}for} each hour $t$ of every representative day $o$. Constraints \eqref{eqn:3d}-\eqref{eqn:4d} enforce permissible frequency response limits pertaining to RoCoF, frequency nadir, and quasi-steady-state frequency deviation, respectively, whereas slack variables $\Delta p^\mathrm{b/s}_{to\psi}$ are used to identify the violations of transient security limits at a specific hour and iteration. Accordingly, \eqref{eqn:1d} provides the minimum change needed in the scheduled power exchange with the main grid from the first stage to ensure frequency security. After solving \eqref{dyn_model} {\color{black}for each considered time step} at each iteration $\psi$, the value of $\Delta p^\mathrm{b/s}_{to\psi}$ is used to modify and tighten the power exchange limits with the main grid at the next iteration ($\psi+1$), as follows: 
\begin{subequations}
\begin{align}
\overline{p}^\mathrm{\: b}_{to(\psi + 1)} = p^\mathrm{b}_{to\psi} - \alpha \Delta p^\mathrm{b}_{to\psi}, \quad \forall t\in\Omega^\mathrm{T},o\in\Omega^\mathrm{O},\label{eqn:minbound}\\
 \overline{p}^\mathrm{\: s}_{to(\psi + 1)} = p^\mathrm{s}_{to\psi} - \alpha \Delta p^\mathrm{s}_{to\psi}, \quad \forall t\in\Omega^\mathrm{T},o\in\Omega^\mathrm{O}.\label{eqn:maxbound}
\end{align}
 \end{subequations}
{\color{black}The scaling factor $\alpha$ is used to apply a less conservative bound modification to account for intertemporal power shifting and investment candidates with frequency support. 
In this work, a value of $\alpha\in[0.5,0.7]$ was adopted.}

\section{Case Study}\label{sec: case studies}
\vspace{-0.1cm}  
\subsection{System Description} \label{subsec: description}
\begin{figure}[t!] 
	\centering
    \includegraphics[clip, trim=3.0cm 10.9cm 7cm 4.7cm,width=0.8\columnwidth]{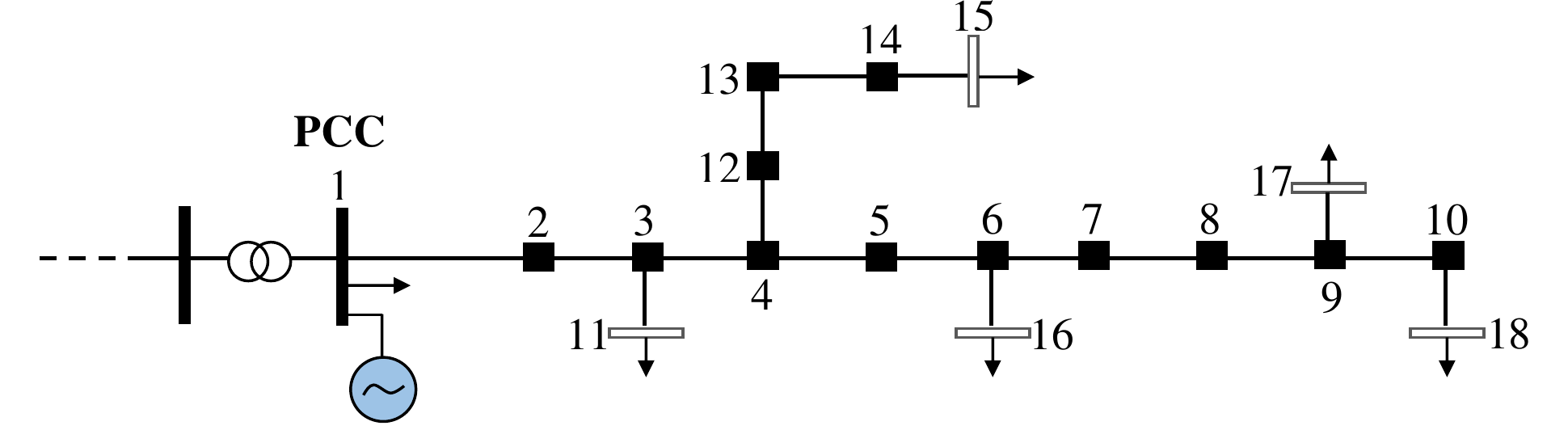}
	\caption{Modified CIGRE European low voltage network.}
	\label{fig:test nw}
	\vspace{-0.35cm}
\end{figure}

A modified CIGRE residential low-voltage network~\cite{cigre}, illustrated in Fig.~\ref{fig:test nw}, is used to analyze the performance of the proposed planning tool. It is assumed that one SG is {\color{black}already preset} at PCC ($\mathrm{SG}_1$) and the investment candidates comprise one SG ($\mathrm{SG}_2$) and three PV CIGs (i.e., $\mathrm{PV}_1$ and $\mathrm{PV}_2$ interfaced via {\it grid-supporting} converters, and $\mathrm{PV}_3$ operating in {\it grid-feeding} mode with fixed power output). The fundamental control parameters {\color{black} obtained from~\cite{UrosUC}} and investment costs ({\color{black} derived from~\cite{6920097}}) of all generators are provided in Table~\ref{tab:gen_param}, while system operation costs are given in Table~\ref{tab:op_cost}. Moreover, $50\,\%$ of nominal load connected at node $1$ is shiftable, whereas high priority critical loads are connected at nodes $15$ and $16$. The patterns of loads and PV generation in Texas during $2016$~\cite{SAM} are used to obtain representative days through $k$-means clustering. Note that all {\color{black}representative days for loads and PV generations as well as} line parameters are provided in the Appendix. The transient security constraints are enforced through thresholds imposed on RoCoF $(\dot{f}_\mathrm{lim}=2\,\mathrm{Hz/s})$, frequency nadir $(\Delta f_\mathrm{lim} =0.8\,\mathrm{Hz})$, and quasi-steady-state frequency deviation $(\Delta f_\mathrm{ss,lim} =0.2\,\mathrm{Hz})$. The implementation was done in \textsc{MATLAB}, with the optimization model formulated in \textsc{YALMIP}~\cite{Lofberg2004} and solved by \textsc{Gurobi}~\cite{gurobi}. 

To analyze the performance of the proposed planning tool, three cases are considered: (1) MG planning without any islanding constraints; (2) MG planning with only static islanding constraints; and (3) MG planning with static and transient frequency islanding constraints.


\begin{table}[!t]
    \renewcommand{\arraystretch}{1.2}
    \centering
    \caption{Generator Control Parameters and Investment Costs}
    \label{tab:gen_param}
    \resizebox{1\linewidth}{!}{%
     \begin{tabular}{c|c|c|c|c|c}
      \toprule
         & \cellcolor{gray!10} $\mathrm{SG}_1$ & \cellcolor{blue!10} $\mathrm{SG}_2$ & \cellcolor{blue!10} $\mathrm{PV}_1$ & \cellcolor{blue!10} $\mathrm{PV}_2$ & \cellcolor{blue!10} $\mathrm{PV}_3$   \\
    \midrule\midrule
        {\color{black}Annualized} investment cost $(\$)$ & \cellcolor{gray!10} - & \cellcolor{blue!10} $40\,000$  & \cellcolor{blue!10} $70\,000$ & \cellcolor{blue!10} $65\,000$ & \cellcolor{blue!10} $60\,000$\\
            Capacity $(\SI{}{\kilo\watt})$  & \cellcolor{gray!10} $280$  & \cellcolor{blue!10} $350$ & \cellcolor{blue!10} $350$ & \cellcolor{blue!10} $350$ & \cellcolor{blue!10} $350$ \\
     Node  & \cellcolor{gray!10} $1$  & \cellcolor{blue!10} $15$ & \cellcolor{blue!10} $17$ & \cellcolor{blue!10} $11$ & \cellcolor{blue!10} $18$ \\
    $M\,(\SI{}{\second})$ & \cellcolor{gray!10} $14$ & \cellcolor{blue!10} $14$ & \cellcolor{blue!10} $14$  & \cellcolor{blue!10} - & \cellcolor{blue!10} -\\
    $D\,(\mathrm{p.u.})$ & \cellcolor{gray!10} $25$ & \cellcolor{blue!10} $25$ & \cellcolor{blue!10} $30$ & \cellcolor{blue!10} - & \cellcolor{blue!10} -\\
    $K\,(\mathrm{p.u.})$ & \cellcolor{gray!10} $1$ & \cellcolor{blue!10} $1$ & \cellcolor{blue!10} $1$  & \cellcolor{blue!10} $1$ & \cellcolor{blue!10} -\\
    $R\,(\mathrm{p.u.})$ & \cellcolor{gray!10} $0.03$ & \cellcolor{blue!10} $0.03$ & \cellcolor{blue!10} - & \cellcolor{blue!10} $0.05$ & \cellcolor{blue!10} -\\
    $F\,(\mathrm{p.u.})$ & \cellcolor{gray!10} $0.35$ & \cellcolor{blue!10} $0.35$ & \cellcolor{blue!10} - & \cellcolor{blue!10} - & \cellcolor{blue!10} -\\
    \bottomrule 
    \rowcolor{gray!10}
    \multicolumn{6}{@{}p{1.5in}}{\footnotesize Existing generator}\\
    \rowcolor{blue!10}
    \multicolumn{6}{@{}p{1.5in}}{\footnotesize Candidate generators}
    \end{tabular}
 }
   \vspace{-0.35cm}
   \end{table}

\begin{table}[!t]
\renewcommand{\arraystretch}{1.2}
 \centering
    \caption{System operation costs}
    \label{tab:op_cost}
\begin{tabular}{c|c|c|c|c|c|c}
\toprule
 \begin{tabular}[c]{@{}c@{}}Import\\ (\$/kWh) \end{tabular} &
 \begin{tabular}[c]{@{}c@{}}Export\\ (\$/kWh) \end{tabular} &
 \begin{tabular}[c]{@{}c@{}}SG\\ (\$/kWh) \end{tabular} &
  \begin{tabular}[c]{@{}c@{}}PV\\ (\$/kWh) \end{tabular} &
  \begin{tabular}[c]{@{}c@{}}Demand\\ shift\\ penalty \\  (\$/kWh) \end{tabular} &
  \multicolumn{2}{l}{\begin{tabular}[c]{@{}c@{}}Demand \\disconnection\\    penalty  \\ (\$/kWh)\end{tabular}} \\\hline
  $30$ &
  $15$ &
  $60$ &
  $0$ &
  $100$ &
  \multicolumn{2}{l}{$(150-200)^*$} \\
  \bottomrule
   \multicolumn{6}{@{}p{3.5in}}{\footnotesize *Based on the level of demand criticality, only in islanded mode}
\end{tabular}
 \vspace{-0.55cm}
\end{table}   


\begin{table}[!b]
    \renewcommand{\arraystretch}{1.2}
    \centering
    \vspace{-0.55cm}
    \caption{Cost comparison with variation in main grid capacity for Case~1{\color{black}: MG Planning without Islanding Constraints}}
    \label{tab:case1}
    \resizebox{1\linewidth}{!}{%
    \begin{tabular}{c|c|c|c|c}
    \toprule
    \begin{tabular}[c]{@{}c@{}}Main grid \\ capacity $(\mathrm{kW})$\end{tabular} & \begin{tabular}[c]{@{}c@{}}Investment costs \\\& decisions $(\$)$\end{tabular} & \begin{tabular}[c]{@{}c@{}}Operational\\costs $(\$)$\end{tabular} & \begin{tabular}[c]{@{}c@{}}Total\\costs (\$)\end{tabular} & \begin{tabular}[c]{@{}c@{}}Installed \\ capacity $(\mathrm{kW})$ \end{tabular}  \\ \hline
    Unlimited & 0 & 56\,394 & 56\,394 & 280   \\ \hline
    250 & 0 & 77\,795 & 77\,795 & 280 \\ \hline
    150 & 60\,000 ($\mathrm{PV}_3$) & 47\,473 & 107\,473 & 630   \\ \hline
    \end{tabular}%
    }
\end{table}


\begin{table*}[!ht]
    \renewcommand{\arraystretch}{1.2}
    \centering
    \caption{Planning Costs for Case~2 (Final Cost in Blue) and Case~3 (Final Cost in Green) {\color{black}and Aggregated Corrective Power Deviations Including} Four Representative Days.}
    \label{tab:cost2_3}
    \begin{tabular}{c|c|c|c|c|c|c|c}
    \toprule
    \begin{tabular}[c]{@{}c@{}}$\psi$\end{tabular} & \begin{tabular}[c]{@{}c@{}}Investment costs\\\& decisions $(\$)$ \end{tabular} & \begin{tabular}[c]{@{}c@{}}Operational \\ costs $(\$)$\end{tabular} & \begin{tabular}[c]{@{}c@{}} {\color{black}Demand}\\{\color{black}shift}\\{\color{black}penalty} $(\$)$\end{tabular} & \begin{tabular}[c]{@{}c@{}} {\color{black}Demand} \\ {\color{black}disconnection} \\ {\color{black}penalty $(\$)$}\end{tabular} & \begin{tabular}[c]{@{}c@{}}Total\\costs (\$)\end{tabular} & \begin{tabular}[c]{@{}c@{}} {\color{black}Import power} \\{\color{black}deviation }\\{\color{black} $(\mathrm{kW})$}\end{tabular} & \begin{tabular}[c]{@{}c@{}} {\color{black}Export power}  \\ {\color{black}deviation }\\ {\color{black}  $(\mathrm{kW})$}\end{tabular}\\ \hline
    1 & 128\,000 ($\mathrm{PV}_2$, $\mathrm{PV}_3$ + Lines 1-2, 2-3, 3-11) & 96\,956 & 3\,613 & 5\,548 & \cellcolor{blue!25} 224\,956 & 2\,902 & 1\,431\\ \hline
    2 & 127\,000 ($\mathrm{PV}_2$, $\mathrm{PV}_3$ + Lines 1-2, 2-3) & 113\,872 & 8\,543 & 5\,337 & 240\,872 & 871 & 429\\ \hline
    3 & 127\,000 ($\mathrm{PV}_2$, $\mathrm{PV}_3$ + Lines 1-2, 2-3)  & 118\,924 & 8\,796 & 5\,081 & 245\,924 & 261 & 129\\ \hline
    4 & 127\,000 ($\mathrm{PV}_2$, $\mathrm{PV}_3$ + Lines 1-2, 2-3)  & 120\,423 & 8\,572 & 5\,334 & 247\,423 & 78 & 39\\ \hline
    5 & 127\,000 ($\mathrm{PV}_2$, $\mathrm{PV}_3$ + Lines 1-2, 2-3)  & 120\,890 & 8\,805 & 5\,081 & \cellcolor{green!25} 247\,890 & 23 & 12\\ \hline
    \end{tabular}%
    \vspace{-0.25cm}
\end{table*}

\vspace{-0.2cm}
\subsection{Cost Analysis}
In this analysis, the costs of {\color{black}the} three aforementioned case studies are compared under the consideration of four representative days. Let us first study Case~1, with the respective costs under different capacity limits of the main feeder listed in Table~\ref{tab:case1}. Understandably, the MG mainly relies on more affordable power provided by the main grid instead of dispatching $\mathrm{SG}_1$ installed at PCC. Under the unlimited import capacity from the main grid, investments {\color{black} in local generation} are not economical due to the low cost of imported power. However, with the introduction of grid capacity limits (e.g., in instances of net load growth and faults experienced in the network), the operational costs increase as a result of {\color{black}the} MG relying on the more expensive $\mathrm{SG}_1$ at PCC. {\color{black}Further} reduction of grid capacity {\color{black}finally} leads to {\color{black}the} installment of $\mathrm{PV}_3$, as it yields higher investment but lower operational costs compared to $\mathrm{SG}_2$ and thus significantly reduces the overall operational costs. 

The variation between investment and operational costs for Cases~2 and 3 is provided in Table~\ref{tab:cost2_3},  where the optimal solution at iteration $\psi=1$ corresponds to the optimal costs of Case~2. The MG requires higher reliability {\color{black}in Case~2 compared to Case~1} in order to minimize the loss of load under static security constraints, whereas in Case~3 the survivability and resilience of the MG are also considered by including the transient security constraints. To ensure the MG resilience, higher investment and operational costs are enforced in both of these case studies compared to Case~1 due to inclusion of static and transient islanding constraints. Indeed, a $400\,\%$ cost increase for Case~2 is observed, with a further $10\,\%$ increase for Case~3. In both of these cases, the installation of renewable DERs reduces the total costs despite the significantly higher underlying investment costs. More precisely, renewable DERs contribute to increased line flows and power export to the main grid, thus necessitating a greater network capacity indicated by the upgrade of {\color{black}the} lines between nodes (1-2) and (2-3). However, in turn, the MG adequacy improves with installing renewable DERs, reflected in the reduction of lost load and ensuring that critical loads are supplied even during emergency islanding situations. 

Focusing on Case~3, it is noticeable that operational costs increase at each iteration due to the use of expensive SGs and flexible loads to mitigate the feasibility violation. 
Nevertheless, at instances where operational flexibility alone cannot guarantee security, more units are installed. Finally, it can be seen that tightening of the power exchange limits (and thus the power export) with the main grid alleviates some of the necessary network investments (e.g., line (3-11) for iterations 1 and 2). 

\vspace{-0.25cm}
\subsection{Transient Security Analysis}

\begin{figure}[t!] 
	\centering
    \scalebox{0.9}{\includegraphics[]{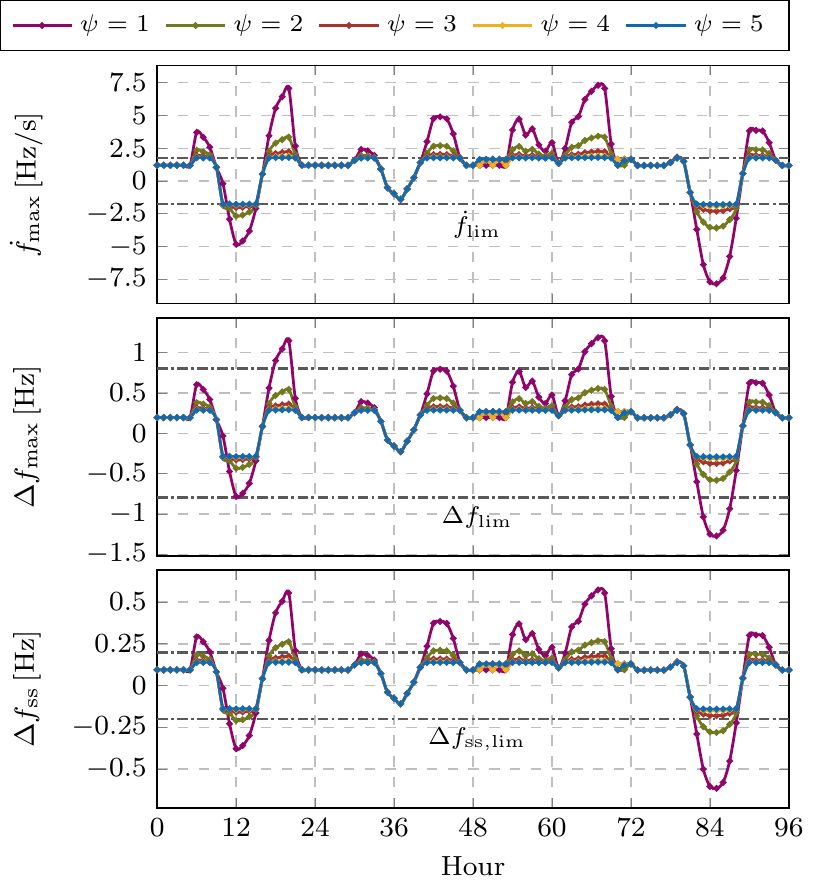}}
    \vspace{-0.25cm}
	\caption{System performance in terms of the transient frequency metrics for different iterations of the proposed algorithm including four representative days.}
	\label{fig:dyn_resp}
	\vspace{-0.5cm}
\end{figure}

In Case~3, the MG survivability is ensured by meeting the prescribed transient security criteria. In the second stage of the algorithm (see Fig.~\ref{fig:methodology}), the slack variable $\Delta p_{to}^\mathrm{b/s}$ is used to indicate the amount of adjustment needed in the scheduled power exchange with the main grid to satisfy the transient frequency requirements. Fig.~\ref{fig:dyn_resp} showcases the metrics describing the dynamic performance of the MG's CoI after islanding at each hour. {\color{black}The} first iteration corresponds to the system response without transient security requirements (Case~2).

A significant improvement is recorded in the maximum RoCoF values, even within a single iteration (e.g., reduction from $8\,\mathrm{Hz/s}$ to $3.5\,\mathrm{Hz/s}$ after the first iteration). Furthermore, each successive iteration reduces the power exchange with the main grid during the hours when security limits are violated until all limits are satisfied. {\color{black} The amount of aggregated corrective power deviations ($\sum_{o\in \Omega^{O}}\sum_{t\in \Omega^{T}}\Delta p_{to}^\mathrm{b/s}$)  in Table~\ref{tab:cost2_3} is monotonically decreased with each iteration until the transient security constraints are fulfilled.}
However, these improvements in terms of security and resilience are achieved at the expense of higher operational costs by dispatching costly SG and flexible loads.
 
It is clear from Fig.~\ref{fig:dyn_resp} that the RoCoF threshold is the most limiting factor for secure transient operation. This is expected, since PV-based CIGs yield a more economic solution but do not provide the same level of inertia as SGs, thus degrading the transient performance. In particular, $\mathrm{SG}_1$ and $\mathrm{PV}_1$ provide both inertia and damping, $\mathrm{PV}_2$ improves damping through droop control, and $\mathrm{PV}_3$ offers no frequency support. Since the inertia and damping contribution of $\mathrm{SG}_1$ and $\mathrm{PV}_2$ do not lead to sufficient transient performance, the reduction in the power exchange with the main grid is needed to ensure a satisfactory response. This is achieved through power provision from $\mathrm{PV}_2$ and $\mathrm{PV}_3$ as well as higher activation of flexible loads.

\vspace{-0.25cm}
\subsection{Sensitivity Analysis}
\subsubsection{Representative Days}

\begin{figure}[t!] 
	\centering
    \scalebox{0.9}{\includegraphics[]{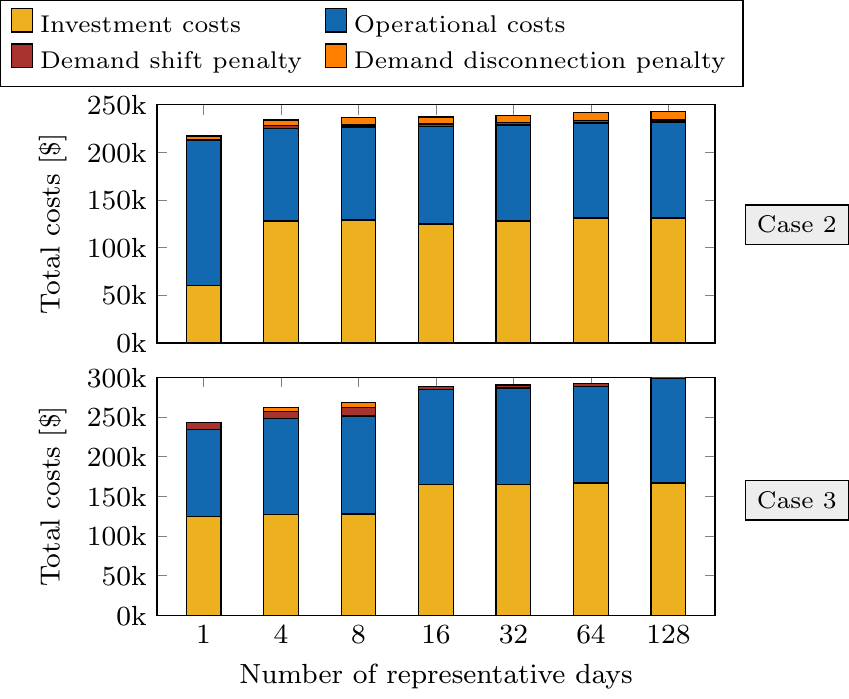}}
    \vspace{-0.25cm}
	\caption{Total costs for different representative days in Cases~2 and 3.}
	\label{fig:rep_days}
	\vspace{-0.45cm}
\end{figure}

The stochasticity of both load and generation profiles affects the planning accuracy, usually resulting in under{\color{black}-} or overestimation. As previously described in Section~\ref{subsec: description}, the load and generation profiles are obtained by utilizing the $k$-means clustering for different representative days. Understandably, the number of considered representative days has a direct impact on the solution of the algorithm. This can be observed in Fig.~\ref{fig:rep_days}, where the total investment and operational costs for Cases~2 and 3 increase with the number of representative days. In particular, employing more representative days provides a better representation of system operation, thus allowing for more accurate estimates of different costs. On the other hand, it also imposes a higher computational burden. In particular, the results in Fig.~\ref{fig:rep_days} indicate that the overall costs plateau for excessive number of representative days, suggesting that the case studies considering up to 16 representative provide a good trade-off between the accuracy of cost estimates and the needed computational effort.
 \vspace{0.3em}
\subsubsection{Operational Flexibility}


\begin{figure}[t!] 
	\centering
    \scalebox{0.9}{\includegraphics[]{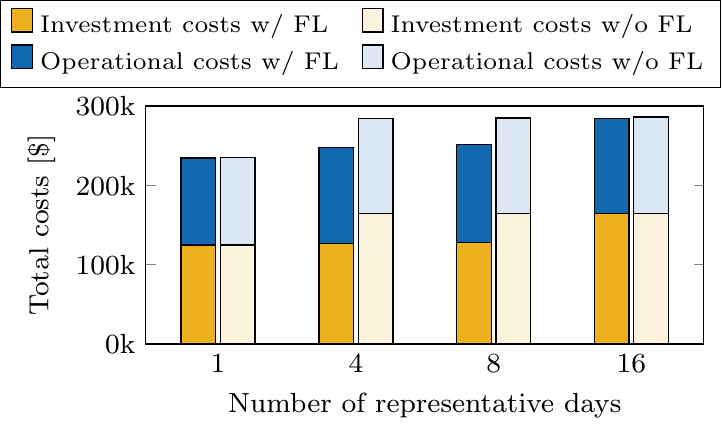}}
    \vspace{-0.25cm}
	\caption{Sensitivity of investment and operational costs to the presence of flexible loads for different representative days in Case~3.}
	\label{fig:flex_sens3}
	\vspace{-0.5cm}
\end{figure}

While flexible loads provide more degrees of freedom for operational planning, they are costly. In spite of their high operational costs, in this study they provide a more affordable option compared to investments in additional generators for improving system flexibility by reducing the peak power exchange with the main grid. Indeed, Table~\ref{tab:cost2_3} shows a successive increase in the use of flexible loads for improving the transient frequency response. This is justified by the fact that flexible loads provide a peak shaving service vital for ensuring survivability during transients.

To this end, Case~3 was studied with and without flexible loads to thoroughly analyze their impact.
In the case of 1 and 16 representative days, the operational costs experience a marginal decrease under the use of flexible loads, whereas the investment costs remain intact, as depicted in Fig.~\ref{fig:flex_sens3}. In contrast, for other representative periods the use of flexible loads leads to lower investment costs, as they alleviate the problems pertaining to adequate power supply. Moreover, in all four cases the total costs increase without the use of flexible loads, thus making their adoption vital for system flexibility and economic operation. The latter aspect is primarily related to the presence of renewable PV units, which allow for the loads to be shifted to periods of higher solar generation. {\color{black}Note that the difference is more prominent in cases with 4 and 8 representative days since the use of flexible loads allows to differ investment decisions}. 

\vspace{-0.25cm}
\section{Conclusion}\label{sec: conclusion}
MGs are expected to play a significant role in increasing the resilience of electric power systems. Their ability to operate in both grid-connected and islanded mode is paramount to their capacity to enhance system reliability. In this paper, the MG investment planning problem under both static and transient islanding constraints is investigated. By explicitly embedding the islanding constraints in the planning problem, the survivability of the system can be guaranteed and the resilience can be assessed as a function of the load supplied in islanded conditions. However, after the islanding event, the transient behavior of the MG is dictated both by the non-linear dynamics and the investment and operation decisions, which poses many challenges concerning the problem formulation. We tackle this problem by proposing an iterative three-stage algorithm that resolves the underlying tractability issues and computational challenges, as well as shows excellent performance on the examined case studies.

Nevertheless, several additional aspects still need to be investigated in ongoing and future work. For instance, the impact of information exchange between different layers of the algorithm on the solution optimality and rate of convergence need to be assessed. Furthermore, a trade-off between the accuracy of the transient response model of the MG and the model complexity should be considered. It is clear though that the need to consider system dynamics within the MG investment and operational decisions is crucial for ensuring system resilience.
\vspace{-0.25cm}

\bibliographystyle{IEEEtran}
\bibliography{bibliography}

\appendices
\section{}  \label{appendixB}

\begin{figure}[h!] 
	\centering
	\scalebox{1.0}{\includegraphics[]{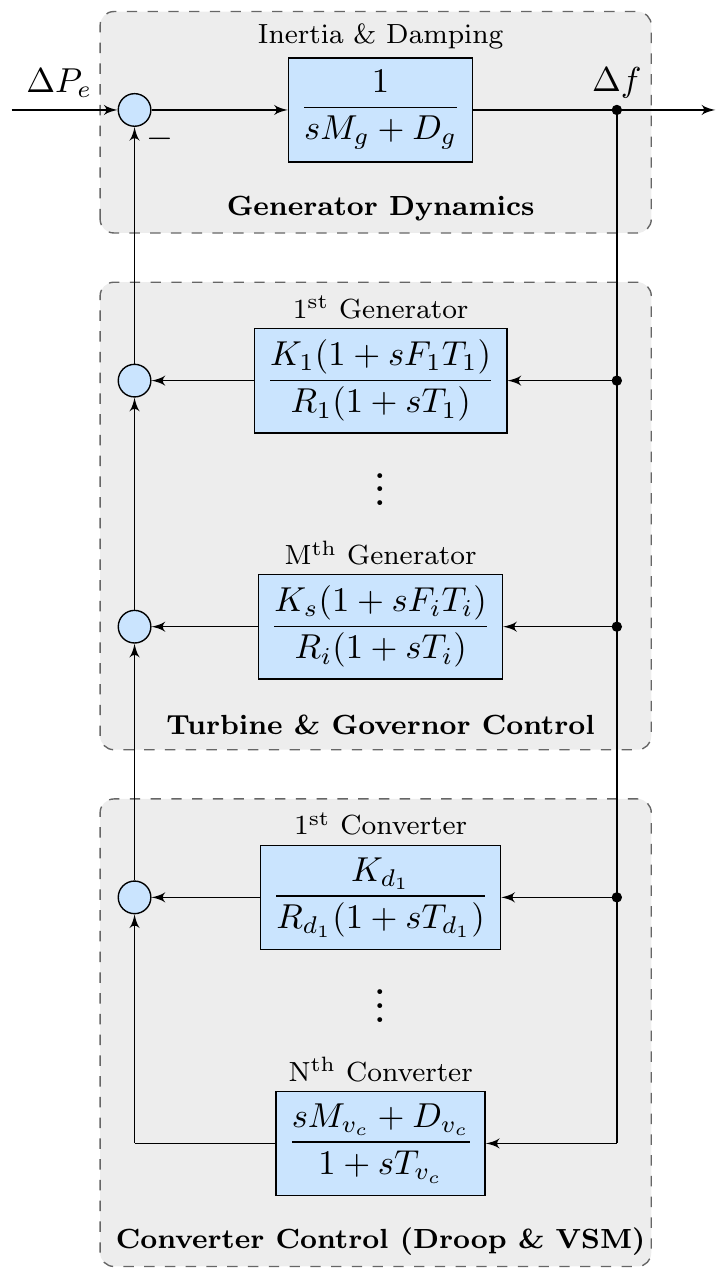}}
	\caption{Uniform system frequency dynamics model.}
	\label{fig:freq_dyn}
\end{figure}

\section{}  \label{appendixA}

\begin{table}[h!]
    \renewcommand{\arraystretch}{1.2}   
    \centering
    \caption{Load Parameters (F: Flexible, C: Constant)}
    \label{load_param}
    \begin{tabular*}{20pc}{c|>{\columncolor{yellow!30}}c|>{\columncolor{yellow!30}}c|>{\columncolor{red!30}}c|>{\columncolor{red!30}}c|>{\columncolor{yellow!30}}c|>{\columncolor{yellow!30}}c}
      \toprule
     Node & 1 & 11 & 15 & 16 & 17 & 18 \\
     Type & F & C & C & C & C & C  \\
     Nominal Load $[\mathrm{kVA}]$ &200 & 15 & 52 & 210 & 35 & 47 \\
     Power factor & 0.95 & 0.95 & 0.95 & 0.85 & 0.95 & 0.95 \\
    \bottomrule
     \rowcolor{red!30}
    \multicolumn{6}{@{}p{1.5in}}{\footnotesize High priority load}\\
    \rowcolor{yellow!30}
    \multicolumn{6}{@{}p{1.5in}}{\footnotesize Low priority load}
    \end{tabular*}
\end{table}

\begin{table}[h]
    \renewcommand{\arraystretch}{1.2}    
    \caption{Line Parameters}
    \label{tab:line param}
    \resizebox{\linewidth}{!}{%
    \begin{tabular}{c|c|c|c|c|c}
    \toprule
    Line & From Node & To Node & Length $[\mathrm{m}]$ & $R\,[\mathrm{p.u.}$ & $X\,[\mathrm{p.u.}$ \\ \midrule\midrule
    1 & 1 & 2 & 35 & 0.010045 & 0.005845 \\ \hline
    2 & 2 & 3 & 35 & 0.010045 & 0.005845 \\ \hline
    3 & 3 & 4 & 35 & 0.010045 & 0.005845 \\ \hline
    4 & 4 & 5 & 35 & 0.010045 & 0.005845 \\ \hline
    5 & 5 & 6 & 35 & 0.010045 & 0.005845 \\ \hline
    6 & 6 & 7 & 35 & 0.010045 & 0.005845 \\ \hline
    7 & 7 & 8 & 35 & 0.010045 & 0.005845 \\ \hline
    8 & 8 & 9 & 35 & 0.010045 & 0.005845 \\ \hline
    9 & 9 & 10 & 35 & 0.010045 & 0.005845 \\ \hline
    10 & 3 & 11 & 30 & 0.03456 & 0.01374 \\ \hline
    11 & 4 & 12 & 35 & 0.04032 & 0.01603 \\ \hline
    12 & 12 & 13 & 35 & 0.04032 & 0.01603 \\ \hline
    13 & 13 & 14 & 35 & 0.04032 & 0.01603 \\ \hline
    14 & 14 & 15 & 30 & 0.03456 & 0.01374 \\ \hline
    15 & 6 & 16 & 90 & 0.1036 & 0.04122 \\ \hline
    16 & 9 & 17 & 30 & 0.03456 & 0.01374 \\ \hline
    17 & 10 & 18 & 30 & 0.03456 & 0.01374 \\ \bottomrule
    \end{tabular}
    }
\end{table}

\newpage
\section{}  \label{appendixC}
\begin{figure}[h!] 
	\centering
	\scalebox{1.0}{\includegraphics[]{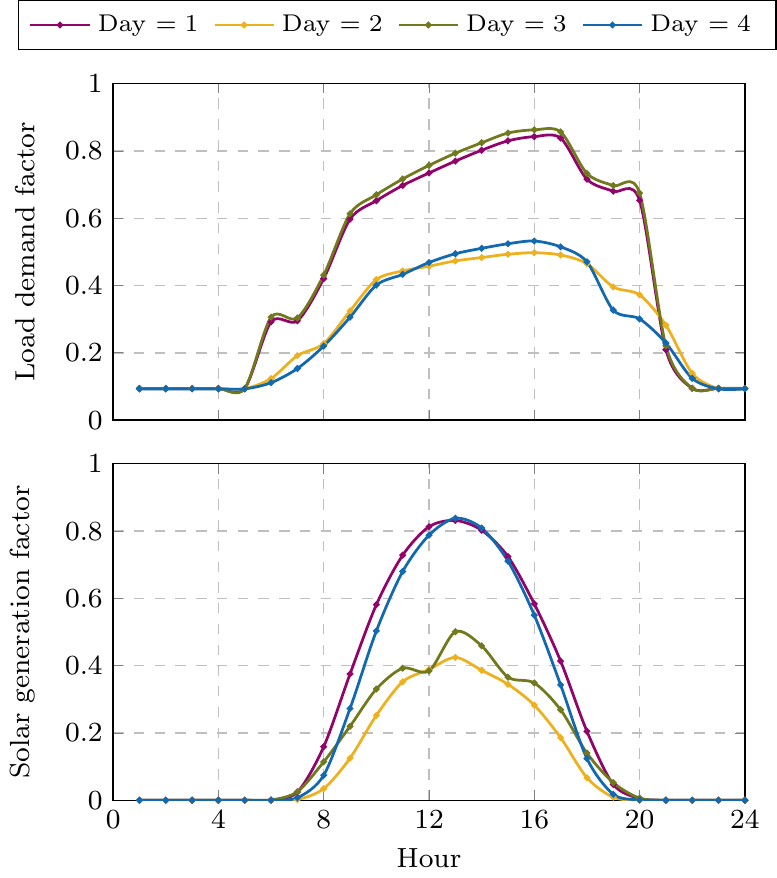}}
	\caption{Demand and solar power generation patterns including four representative days.}
	\label{fig:profiles}
\end{figure}

\end{document}